\def\versiondate{17 Oct. 2016}
\input math10.macros
\magnification=1200
\hsize=6.9truein
\hoffset=-.3truein  %

\let\nobibtex = t
\let\noarrow = t
\input eplain.tex
\beginpackages
\usepackage{url}
\endpackages

\input lrlEPSfig.macros

\checkdefinedreferencetrue
\continuousfigurenumberingtrue
\theoremcountingtrue
\sectionnumberstrue
\forwardreferencetrue
\citationgenerationtrue
\nobracketcittrue
\hyperstrue
\initialeqmacro

\input\jobname.key
\bibsty{../../texstuff/myapalike}

\def\cp{\tau}   %
\def\verts{{V}}
\def\vertex{{V}}
\def\edges{{E}}
\def\tr{\mathop{\rm tr}}  %
\def\Tr{\mathop{\rm Tr}}  %

\def\bfz{{\bf 0}}

\def\bigpip#1{\bigl(\changecomma #1\bigr)}

\def\changecomma#1,{#1,\,}
\def\bigchangecomma#1,{#1,\;}
\def\leftchangecomma#1,{#1,\ }

\def\bfz{{\bf 0}}

\def\ent{{\bf H}}  %
\def\ER(#1,#2){{\cal R}(#1 \leftrightarrow #2)}   %
\def\dom{\succcurlyeq}
\def\LL{{\cal L}}  %
\def\etail#1{#1^-}   %
\def\ehead#1{#1^+}   %

\def\ev#1{{\scr #1}}

\def\BLPSusf{\ref b.BLPS:usf/, hereinafter referred to as BLPS (2001)%
\def\BLPSusf{BLPS \htmllocref{\bibcode{BLPS:usf}}{(2001)}}}

\def\BLPSgip{\ref b.BLPS:gip/, hereinafter referred to as BLPS (1999)%
\def\BLPSgip{BLPS \htmllocref{\bibcode{BLPS:gip}}{(1999)}}}

\ifproofmode \relax \else\head{To appear in {\it Combin. Proba. Comput.}}
{Version of \versiondate}\fi 
\vglue20pt

\title{Comparing Graphs of Different Sizes}

\author{Russell Lyons}

\abstract{We consider two notions describing how one finite graph may be
larger than another. Using them, we prove several theorems for such pairs
that compare the number of spanning trees, the return probabilities of
random walks, and the number of independent sets, among other combinatorial
quantities. Our methods involve inequalities for determinants, for traces
of functions of operators, and for entropy.
}

\bottomII{Primary 
 05C05, %
60C05. %
Secondary
 05C80, %
 05C81. %
}
{Spanning trees, random walks, independent sets, matchings.}
{Research partially supported by NSF grants DMS-0406017 and DMS-1007244,
and Microsoft Research.}

\bsection{Introduction}{s.intro}

How does one compare different graphs? 
If they have the same vertex sets but one edge set contains the other,
then clearly we can say that one graph is larger than the other.
This can lead to inequalities, usually trivial, for combinatorial
quantities associated with the graphs.
But if the numbers of vertices differ, then one might wish to compare those
combinatorial quantities with normalizations that depend on the numbers of
vertices.
In such a case, however,
if one graph can be embedded in another, we are unlikely to
have any such comparison of normalized combinatorial quantities.
Instead, we should demand some sort of uniformity in how one graph can be
embedded in the other.
With appropriate hypotheses of uniformity,
we have found inequalities for
counting any of the following: spanning trees; independent sets; proper
colorings; acyclic orientations; forests; and matchings.
We also have inequalities for random walks and the spectra of the graphs.
However, many questions remain open.

We now describe what we mean by uniformity of embedding.
Let $G$ and $H$ be finite connected
(multi)graphs. We will use $G$ for the larger graph.
The simplest kind of uniformity is
that $H$ \dfn{tiles} $G$, meaning that $G$ contains a
collection of copies of $H$ that cover each vertex of $G$ exactly once.
See \ref f.sqtiles/ for an example.
The case where $H$ tiles $G$
is hardly different from $H$ being a subgraph of $G$ with the same number
of vertices, and will not be discussed further here.

\efiginslabel sqtiles {The $4 \times 4$ portion of the square grid is tiled
by 4 copies of a 4-cycle.} x2 

In general, we define a \dfn{copy} of $H$ to be a subgraph of $G$ that is
isomorphic to $H$.
We will also refer to a copy of $H$ as an embedding of $H$.

The next simplest kind of uniformity is that $G$ has a \dfn{fractional
tiling} by $H$.
This means that there is an integer number of copies of $H$ in $G$ such
that each vertex of $G$ is covered the same number of times by
these copies of $H$. 
An example is in \ref f.fractionaltile/.
This is already nontrivial and will be a common hypothesis in our paper.

\twoefiginsonelabel K4 {y1.5} 3.5 triangle {y1.5} 3 fractionaltile
{The graph $K_4$ is fractionally tiled by $K_3$.} 

Finally, the most general case we will consider is
that $G$ \dfn{dominates} $H$,
written $G \succcurlyeq H$, 
meaning that there is a probability measure on pairs $(X, Y) \in
\vertex(G) \times \vertex(H)$ such that almost surely
there is a rooted isomorphism from 
$(H, Y)$ to a subgraph of $(G, X)$ and such that the marginal distributions
of $X$ and $Y$ are each uniform.
Here, a rooted graph is a pair $(G, o)$ with $o \in \verts(G)$ and a rooted
isomorphism is an isomorphism that carries one root to the other.
The way to think of domination is that $G$ looks bigger than $H$ from the
point of view of a typical vertex.
For some illustrative examples, see Figures 
\briefref f.domnotfrac/ and \briefref f.domsails/.

\midinsert   
      \line{%
      \hbox to 4.25truein{\hfill\Size y1.5 \epsfbox{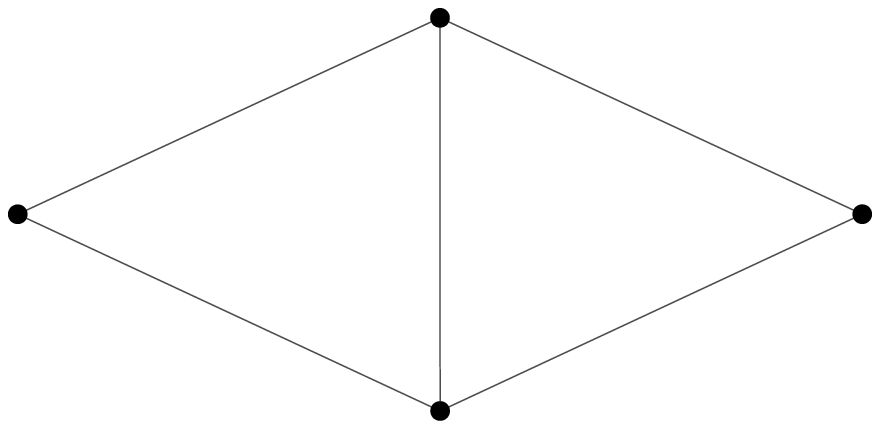}\hfill}%
      \hbox to 2.25truein{\hfill\Size y1.5 \epsfbox{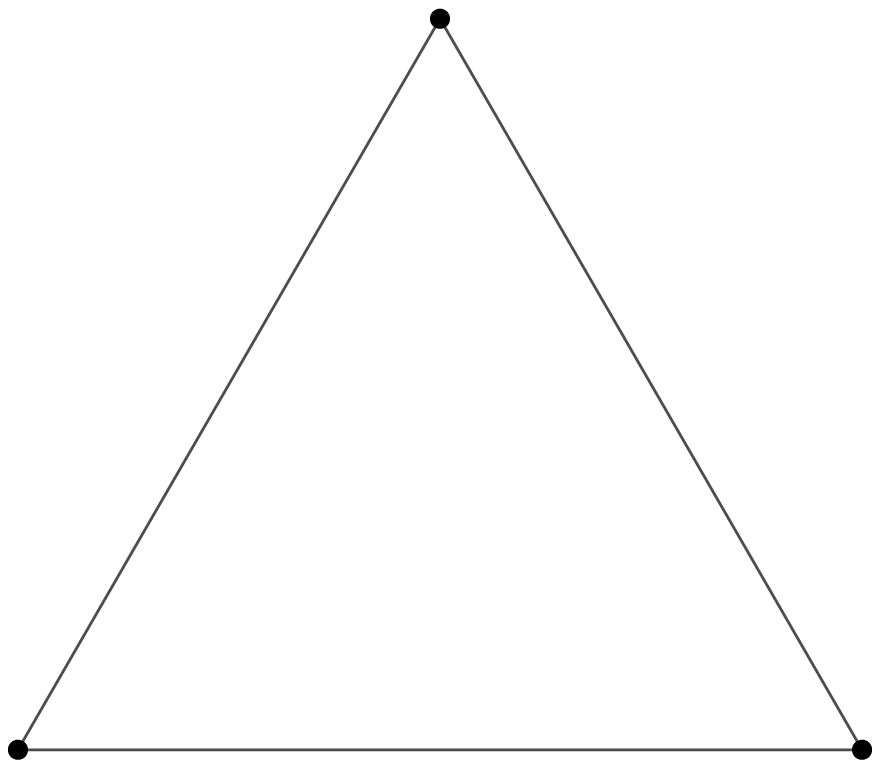}\hfill}}%
\caption{\hsize=3.8truein%
\vtop{\noindent \figlabel{domnotfrac}\enspace 
The graph $G$ on the left dominates the graph $H$
on the right, but $H$ does not fractionally tile $G$.
}}
\endinsert

\twoefiginsonelabel domsails {y1.2} 3.5 sails2 {y1.2} 3 domsails
{The graph on the left dominates the graph on the right.}

We say that a graph is \dfn{transitive} if for every pair of vertices,
there is an automorphism of the graph that takes one to the other.
If $H$ is transitive, then $G \dom H$ iff every vertex of $G$
belongs to a copy of $H$. If $G$ is transitive, then $G \dom H$ iff $G$
contains a copy of $H$. 
In both cases, the independent coupling of roots works.
It is clear that if $H$ fractionally tiles $G$, then $G \dom H$.
Conversely, if $G$ is transitive and dominates $H$, then $H$ fractionally
tiles $G$.

Consider the case where $H$ is a single edge.
Then to say that $G$ dominates $H$ is to say
that $G$ has no isolated vertices; since $G$ is connected, this means that
$G$ has at least two vertices. On the other hand, to say that $H$
fractionally tiles $G$ is to say that there is a spanning ``subgraph" of
$G$ that is regular of
degree at least 1; the reason for the quotes is that the subgraph may
need to use edges of $G$ multiple times and thus be a multigraph even if
$G$ is a simple graph.
See \ref f.domfracedge/ for a comparison.

\midinsert   
      \line{%
      \raise.5truein\hbox to 1.7truein{\hfill\Size x1.5 \epsfbox{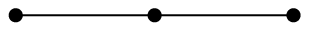}\hfill}%
      \hbox to 2.8truein{\hfill\Size y1.2 \epsfbox{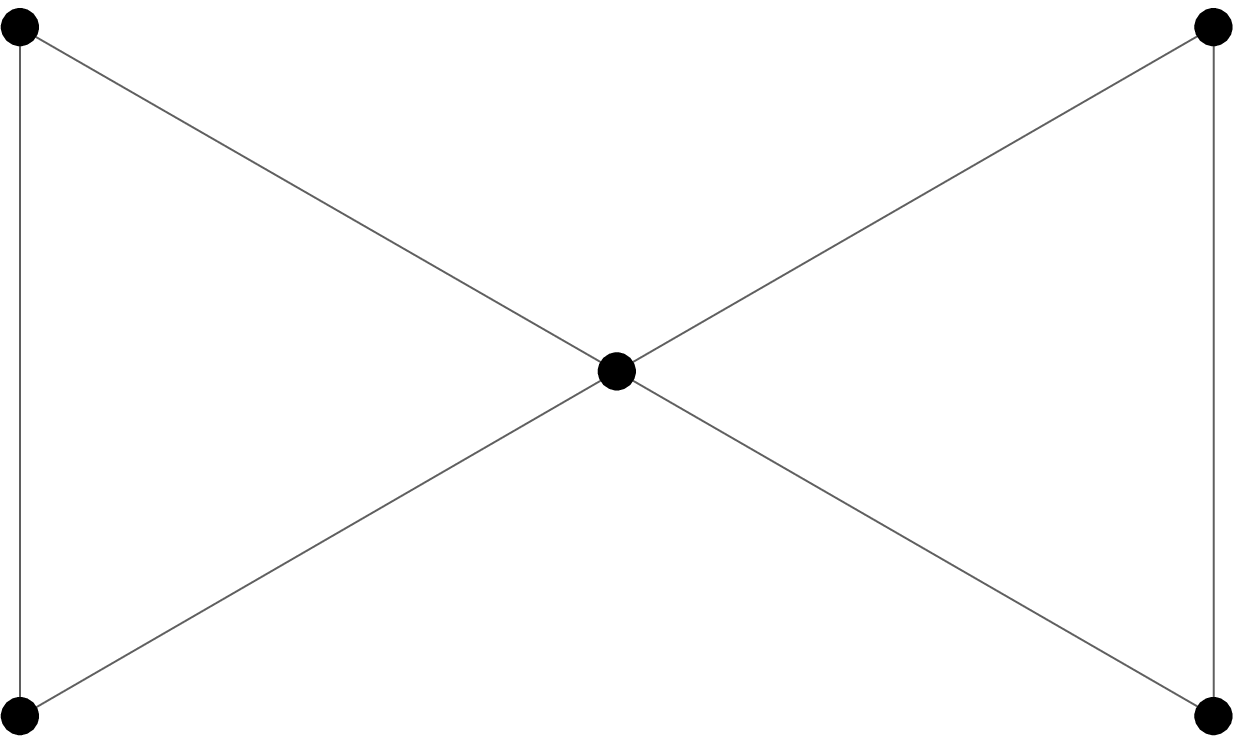}\hfill}%
      \hbox to 2.truein{\hfill\Size y1.2 \epsfbox{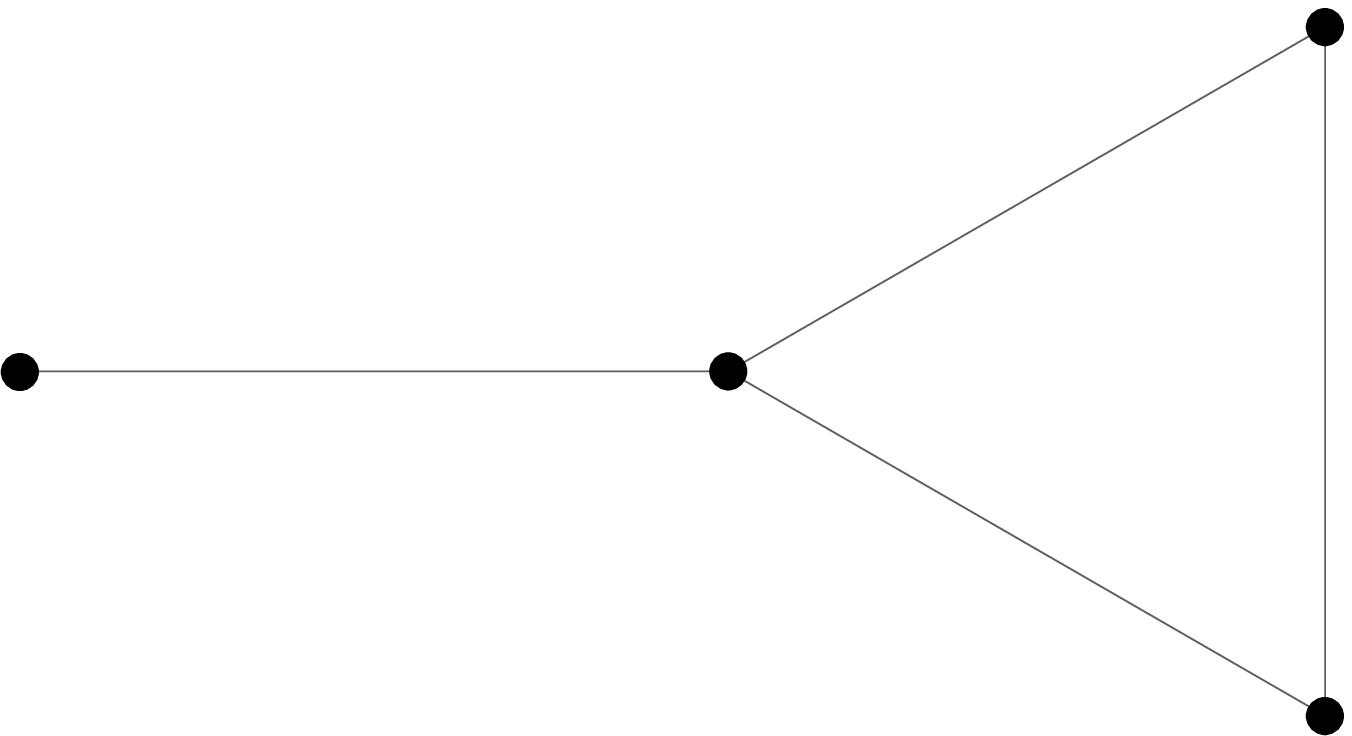}\hfill}}%
\caption{\hsize=5truein%
\vtop{\noindent \figlabel{domfracedge}\enspace 
The graph on the left dominates an edge; an edge fractionally tiles
the graph in the middle and tiles the graph on the right.
}}
\endinsert

In probabilistic language, $G \dom H$ has a simple expression, though we
will not use it. Note that
the set of (isomorphism classes of)
rooted graphs is partially ordered by rooted embedding. This
partial order defines a corresponding notion of stochastic ordering $\dom$
on the set of probability measures on rooted finite graphs.
Let $U(G)$ denote the probability measure on $G$ with a uniformly random
root. Then $G \dom H$ iff $U(G) \dom U(H)$.

Our theme is the following.
Suppose we know an inequality of the form $f(H) \le f(G)$ when $H$ is a
spanning subgraph of $G$. Does it extend with appropriate normalization to
the setting of domination or fractional tiling?

One might consider weighted graphs as well, but in most cases, we will not.
One can also consider random graphs; see Sec.~6 of
\ref b.LPS:GD/ for several such questions. \ref b.LWbuild/ and \ref
b.Janson/ contain further information for random trees.

If $G$ is a connected
graph, we call a subgraph of $G$ a \dfn{spanning tree} if it is
maximal without cycles.
The number of spanning trees of $G$ is denoted $\cp(G)$.

We conjecture the following:

\procl g.main
If $G \dom H$, then 
$$
\cp(G)^{1/|G|} \ge \cp(H)^{1/|H|}
\,.
\label e.spnineq
$$
\endprocl

The infinitary analogue of \ref g.main/
for unimodular probability measures is true, though we will not use it
explicitly; see \ref b.Lyons:trent/.
We will establish several special cases of \ref g.main/.
These proofs will use the Hadamard--Fischer--Koteljanskii inequality for
determinants, which we review in \ref s.back/.
One could extend our considerations to stochastic domination of probability
measures of the form $U(G)$ and $U(H)$ for $G$ and $H$ themselves random,
but for simplicity, we usually avoid such.
The infinitary analogue of \ref g.main/ also
implies that \ref e.spnineq/ holds for certain pairs of large graphs: see
\ref p.large/ below.

After treating spanning trees, we will give inequalities for return
probabilities of continuous-time random walks and for 
eigenvalues when one graph fractionally tiles another.
The main tool here will be a trace inequality for functions of operators.

Our last section presents some easy consequences of Shearer's inequality
about entropy
for a variety of combinatorial quantities, such as counting the number of
independent sets.
The last two sections contain several open questions.

\bsection{Determinant Inequalities}{s.back}

For a square matrix $M$ and a subset $A$ of the indices of its rows and
columns, let $M(A)$ denote the minor of $M$ corresponding to the rows and
columns indexed by $A$.
We use the convention $M(\emptyset) := 1$.
The Hadamard--Fischer--Koteljanskii inequality says that if $M$
is a positive semidefinite matrix, then
$$
M(A) M(B) \ge M(A \cup B) M(A \cap B)
\,.
$$
In other words, $M(\, \cbuldot \,)$ is log-submodular.
It follows that if the index of
each row belongs to precisely $m$ sets $A_i$ (each of arbitrary size), then
$$ 
\prod_i M(A_i) \ge (\det M)^m
\,.
\label e.FHK
$$
Indeed, we simply aggregate repeatedly any pairs of subsets where neither is
contained in the other. Each time we get a larger subset and this can
only end when we have $m$ copies of the entire index set (together with
irrelevant empty sets).

From now on, \dfn{all graphs we consider will be connected without mention}.

One well-known example of a positive semidefinite matrix associated to a
finite graph, $G$, is its \dfn{Laplacian}, $\Delta_G$, whose off-diagonal
entries $\Delta_G(x, y)$ are negative 
the numbers of edges joining $x$ and $y$ and
whose row sums each vanish.
For $W \subset \verts(G)$, denote by $G/W$ the graph obtained from $G$ by
identifying all vertices in $W$ to a single vertex.
By the Matrix-Tree theorem, for any non-empty subset $W \subset
\vertex(G)$, we have 
$$
\cp(G/W) = \Delta_G\big(\vertex(G) \setminus W\big)
\,.
\label e.MTT
$$
Thus, if we denote $G/(\vertex(G) \setminus A)$ by $G_A$, then we have 
$$
\cp(G_A) \cp(G_B) \ge \cp(G_{A \cup B}) \cp(G_{A \cap B})
\label e.step
$$
when $A \cup B$ is a proper subset of $\vertex(G)$.

Inequality \ref e.step/ does not hold when $A \cup B = \vertex(G)$.
For example, take $G$ to be
a path on 3 vertices, $x$, $y$, and $z$, with $y$ the middle vertex.
Let $A := \{x, y\}$ and $B := \{y, z\}$.
Then the left-hand side of \ref e.step/ is 1 while the right-hand side is
2.
However, \ref e.step/ does hold when $A \cup B = \verts(G)$ and
there is an edge between $A \setminus B$ and $B \setminus A$.
Indeed,
such an edge may be subdivided to create a new vertex $x$ that does not
belong to either $A$ or $B$. 
Let the new graph be $G'$ with vertex set $A \cup B \cup \{x\}$.
Note that $\cp(G'_A) = \cp(G_A)$, $\> \cp(G'_B) = \cp(G_B)$, $\> \cp(G'_{A
\cup B}) = 2 > 1 = \cp(G_{A \cup B})$, and $\cp(G'_{A \cap B}) \ge \cp(G_{A
\cap B})$. Thus,
if we apply \ref e.step/ to $G'$ with these same $A$ and $B$, we obtain an
inequality that is stronger than \ref e.step/ applied to $G$, as desired.

As Jeff Kahn noted (personal communication), this extension of \ref e.step/
implies the following inequality: if $A_i \subset \vertex(G)$ and each $x
\in \vertex(G)$ belongs to exactly $m$ sets $A_i$, then 
$$
\prod_i \cp(G_{A_i}) \ge \cp(G)^m
\,,
\label e.stepjk
$$
since if
$A_i \ne \vertex(G)$, then there is an $x \in A_i$ adjacent to some $y
\notin A_i$, whence there is some $A_j$ containing $y$ but not $x$, as $x$
is covered the same number of times as $y$ is. Hence we can apply \ref
e.step/ or its extension repeatedly.

Some final notation for a (possibly disconnected) subgraph $H$ of $G$:
Write $G_H$ for $G_{\vertex(H)}$.
Write $G//H := G/\edges(H)$, the graph obtained from $G$ by contracting
all edges in $\edges(H)$. When $H$ is connected, $G//H =
G/\vertex(H)$.
By $|G|$, we mean $|\vertex(G)|$.

\bsection{Spanning Trees}{s.transt}

In this section, we prove \ref e.spnineq/ when $G \dom H$ and
either $G$ or $H$ is
transitive, and in general when $G$ is fractionally tiled by $H$.

\procl l.easy
For any $G$ and $H \subseteq G$, we have 
$$
\cp(H) \cp(G//H) \le \cp(G)
\,.
$$
\endprocl

\proof
The union of the edges of a spanning tree of $H$ and a spanning tree of
$G//H$ is a spanning tree of $G$. This map from the pairs of spanning trees
of $H$ and the spanning trees of $G//H$ to the spanning trees of $G$
is obviously injective.
\Qed

\procl l.hft
If $G$ is transitive and $H \subseteq G$, then 
$$
\cp(G_H)
\ge
\cp(G)^{|H|/|G|}
\,.
$$
\endprocl

\proof
We prove this by induction on $|H|$.
We may assume that $\emptyset \ne H \ne G$.
Call an image of $H$ under an automorphism of $G$ a \dfn{clone} of $H$.
We claim that there is an edge $e$ such
that we can cover $\vertex(G)$ by clones of $H$, none of which
uses both endpoints of the edge $e$.  Indeed, cover by clones as much of
$\vertex(G)$ as possible without covering all of $\vertex(G)$. Let $o$ be a
vertex that is not covered. Any clone of $H$ that covers $o$ covers all
other uncovered vertices of $G$.
Further, some vertex $x \in \verts(H)$ has degree less than its degree in
$G$. Choose an automorphism of
$G$ that
maps $x$ to $o$ and use the corresponding clone of $H$ to
finish the cover of $\vertex(G)$. Let $o'$ be a neighbor of $o$ that is not
contained in this last clone of $H$.
The edge $e = (o, o')$ is the edge we desire.
Let $k$ be the total number of clones $H_i$ of $H$ used in this cover of
$\vertex(G)$.

Now let $G'$ be $G$ subdivided at $e$ by a new vertex, $z$. We have
$\vertex(G' \setminus z) = \bigcup_{i=1}^k \vertex(H_i)$.
Note that for each $i$, $\>G'_{H_i}$ is isomorphic to $G_{H_i}$, possibly
plus a loop,
because $H_i$ does not
contain both $o$ and $o'$.
Furthermore, $G_{H_i}$ is isomorphic to $G_H$ by definition of clone.
We may assume that each $H_i$ has a vertex not belonging to
$L_i := \bigcup_{j < i} H_j$.
By \ref e.step/, we have 
$$
\cp(G'_{L_{i}}) \cp(G'_{H_i}) 
\ge
\cp(G'_{L_{i+1}}) \cp(G'_{K_i})
\,,
$$
where $K_i := L_{i} \cap H_i \subsetneq H_i$.
Since $\vertex(L_{k+1}) = \vertex(G' \setminus z)$, we have $G'_{L_{k+1}} = G'$.
In addition, since $K_i \subset H_i$, we have $G'_{K_i}$ is $G_{K_i}$,
possibly plus
a loop.
Thus, multiplying together the above inequalities for $1 \le i \le k$ and
cancelling common terms on both sides yields
$$
\cp(G_H)^k
=
\prod_{i \ge 1} \cp(G'_{H_i})  
\ge
\cp(G') \prod_{i > 1} \cp(G'_{K_i})
=
\cp(G') \prod_{i > 1} \cp(G_{K_i})
\ge
\cp(G) \prod_{i > 1} \cp(G_{K_i})
\,.
$$

Since $|K_i| < |H_i|$,
the inductive hypothesis gives $\cp(G_{K_i}) \ge
\cp(G)^{|K_i|/|G|}$.
Now $k |H| = |G| + \sum_{i > 1} |K_i|$, whence 
$$
\cp(G_H)^k 
\ge
\cp(G)^{k |H|/|G|}
\,,
$$
which is the desired inequality.
\Qed

Strict inequality holds in \ref l.hft/ when $|G| > |H| \ge 1$ and
$G$ contains no cut-edge, since
in that case, $\cp(G') > \cp(G)$ in the proof.

\procl t.Gtrns
If $G$ is transitive, then \ref e.spnineq/ holds with strict inequality when
$G$ contains no cut-edge and $G \ne H$.
\endprocl

\proof
Let $K$ be the complement of the vertices in a copy of $H$ in $G$.
The previous two lemmas give
$$
\cp(G)^{|K|/|G|} \le \cp(G_K)
=
\cp(G//H)
\le
{\cp(G) \over \cp(H)}
\,.
$$
Since $|K| + |H| = |G|$, the desired inequality follows.
\Qed

We now prove that \ref g.main/ holds when $H$ is transitive.

\procl t.Htrns
If $H$ is transitive, then \ref e.spnineq/ holds.
\endprocl

\proof
We prove this by induction on $|G|$.
If $|G| = |H|$, then $G$ contains a copy of $H$ and the result is trivial.
Otherwise,
cover $\vertex(G)$ by copies $H_i$ of $H$ for $1 \le i \le k$.
For $1 < i < k$, we may assume that $H_i$ has a vertex either in or
adjacent to $H_j$ for some $j < i$.
Let $V' := \bigcup_{1 \le i < k} \vertex(H_i)$ and
$G'$ be the graph spanned by $V'$.
Then $G'$ is connected.
We may also assume that $H' := H_k$ has a vertex not in $V'$.
Since $H$ is transitive and $G'$ is covered by copies of $H$, we know that
$G' \dom H$, whence our inductive hypothesis says that 
$$
\cp(G')^{1/|G'|} \ge \cp(H)^{1/|H|}
\,.
$$
If $H'$ does not contain a vertex in $V'$, then
$$
\cp(G)
\ge
\cp(G') \cp(H')
\ge
\cp(H)^{|G'|/|H|} \cp(H')
=
\cp(H)^{|G|/|H|}
\,.
$$
If $H'$ does contain a vertex in $V'$, then
note that $G//G' = G/V'$ is isomorphic, up to loops,
to $H'/(V' \cap H') = H'_{H' \setminus V'}$.
Thus, the previous two lemmas give
$$\eqaln{
\cp(G)
&\ge
\cp(G') \cp(G//G')
\ge
\cp(H)^{|G'|/|H|} \cp(H'_{H' \setminus V'})
\cr&\ge
\cp(H)^{|G'|/|H|} \cp(H)^{|H' \setminus V'|/|H|}
=
\cp(H)^{|G|/|H|}
\,.
}$$
Both cases together complete the induction.
\Qed

Note that the same proof shows that if $H$ is any 
graph such that every (connected)
subgraph $K \subset H$ satisfies $\cp(H_K) \ge \cp(H)^{|K|/|H|}$, then for
every $G$ each of whose vertices belongs to a copy of $H$ [in particular,
if $G \succcurlyeq H$], we
have $\cp(G)^{1/|G|} \ge \cp(H)^{1/|H|}$. Many small non-transitive
graphs $H$ can be shown to have this property.

We also show that pairs of large graphs tend to satisfy \ref g.main/.
Write $\|\,\cbuldot\,\|$ for the usual total-variation norm of signed
measures.
Also, write $U_r(G)$ for the distribution of (the isomorphism class of)
the rooted ball of radius $r$
about a uniform random root of $G$.

\procl p.large
Suppose that $D, r < \infty$ and $\epsilon> 0$. There is some $k < \infty$
such that if $G \dom H$, $\>|H| \ge k$, all degrees in $G$ are at most $D$,
and $\|U_r(G) - U_r(H)\| \ge \epsilon$, then $\cp(G)^{1/|G|} \ge
\cp(H)^{1/|H|}$.
\endprocl

\proof
If not, then there is a sequence $G_n \dom H_n$ with $|H_n| \to\infty$, all
degrees of $G_n$ are at most $D$, $\>\|U_r(G_n) - U_r(H_n)\| \ge \epsilon$,
and $\cp(G_n)^{1/|G_n|} < \cp(H_n)^{1/|H_n|}$.
By compactness, there is a subsequence, which for simplicity of notation we
take to be the whole sequence, such that $U(G_n)$ weakly converges to
some probability measure $\mu$ on rooted graphs and $U(H_n)$ weakly
converges to some probability measure $\nu$ on rooted graphs.
Then $\|\mu - \nu\| \ge \epsilon$.
By Theorem 3.2 of \ref b.Lyons:est/ and
Theorem 3.3 of \ref b.Lyons:trent/, we have that $\lim_{n \to\infty}
\cp(G_n)^{1/|G_n|} > \lim_{n \to\infty} \cp(H_n)^{1/|H_n|}$, a
contradiction.
\Qed

We remark that weaker assumptions suffice in place of the bounded degree
assumption; as long as tightness and bounded average log degree
hold for a class of graphs, the same
argument works. See Section 3 of \ref b.BLS:urt/ for a discussion of
tightness.

We owe the following result to Jeff Kahn.

\procl t.spnfrac
If $H$ fractionally tiles $G$, then \ref e.spnineq/ holds.
\endprocl

\proof
Let $H_i$ be the copies of $H$ that fractionally tile $G$ and $A_i :=
\vertex(G) \setminus \vertex(H_i)$ for $1 \le i \le m$.
We have by \ref l.easy/ that 
$$
\cp(G)^m \ge \prod_i \cp(H_i) \cp(G//H_i)
\,.
$$
Since $H_i$ is connected, we have $G//H_i = G_{A_i}$, so now we can apply
\ref e.stepjk/.
Each vertex of $G$ appears $m\big(|G|-|H|\big)/|G|$ times in some $A_i$, whence
$$
\cp(G)^m \ge \cp(G)^{m\onebig(|G|-|H|\onebig)/|G|} \prod_i \cp(H_i) 
=
\cp(H)^m \cp(G)^{m\onebig(|G|-|H|\onebig)/|G|}
\,.
$$
This gives the desired inequality.
\Qed

In summary, we have proved that \ref g.main/ holds under any of the
following additional hypotheses:
if either $G$ or $H$ is transitive;
if $H$ fractionally tiles $G$;
if $H$ is any graph such that every (connected)
subgraph $K \subset H$ satisfies $\cp(H_K) \ge \cp(H)^{|K|/|H|}$;
if $H$ is sufficiently large and $G$ and $H$ are sufficiently distinct (see
\ref p.large/ for details).

\bsection{Fractional Tiling and Random Walks}{s.ftr}

For a continuous-time random walk on a weighted simple graph $G$, let
$p_t(x; G)$ denote the probability that a random walk started at $x$ is at
$x$ at time $t$.
If $\Delta_G$ is the corresponding Laplacian, i.e., $\Delta_G(x, y) :=
- w(e)$ when $e$ is an edge joining $x$ and $y$ with weight $w(e)$,
all other off-diagonal
elements of $\Delta_G$ are 0, and the row sums are 0, then $p_t(x; G)$ is
the $(x, x)$-entry of $e^{-t \Delta_G}$.

We would like to prove that if $G$ dominates $H$, then
for all $t > 0$,
$$
{1 \over |G|} \sum_{x \in \vertex(G)} p_t(x; G)
\le
{1 \over |H|} \sum_{x \in \vertex(H)} p_t(x; H)
\,.
\label e.ftr0
$$
It is easy to see that this inequality holds near 0 and near $\infty$.
One motivation is the following open problem of Fontes and Mathieu
(personal communication). 
Suppose that $G$ is a fixed Cayley graph and $w_1$, $w_2$ are two random
fields of
positive weights on its edges with the following properties: 
Each field $w_i$ has an invariant law and a.s.\ $w_1(e) \ge w_2(e)$
for each edge $e$.
Does it follow that $\Ebig{p_{1, t}(o; G)} \le \Ebig{p_{2, t}(o, G)}$ for all $t
> 0$, where $p_{i, t}$ denotes the return probability to a fixed vertex $o$
at time $t$ with the weights $w_i$?
This is known to be true for amenable $G$ (\ref B.FontesMathieu/) and also
when the pair $(w_1, w_2)$ is invariant (\ref B.AL:urn/).
The problems for finite graphs and for infinite Cayley graphs 
are quite similar in that both try to compare different normalized traces.

We prove a partial result, namely, that \ref e.ftr0/ holds when $H$
fractionally tiles $G$. 

\procl t.ftr0
If $G$ is fractionally tiled by $H$, then for continuous-time simple
random walk, we have
for all $t > 0$,
$$
{1 \over |G|} \sum_{x \in \vertex(G)} p_t(x; G)
\le
{1 \over |H|} \sum_{x \in \vertex(H)} p_t(x; H)
\,.
$$
Equality holds iff $G = H$.
\endprocl

In fact, a somewhat weaker condition suffices: a number of different graphs
can be used inside $G$ as long as their average is at most $G$ in a certain
sense, as we formalize next.
The equality condition of \ref t.ftr0/
arises from the proof of \ref t.ftr/:
we have strict inequality in
\ref e.bigger/ if $G \ne H$.
In the following result, the case when $k = m = 1$ is due to Benjamini
and Schramm; see Theorem 3.1 of \ref b.HeicklenHoffman/.

\procl t.ftr
Let $G$ be a graph with positive weights $w$ on its edges.
Suppose that $H_i$ is a subgraph of $G$ with positive weights $w_i$ on its
edges for $i = 1, \ldots, k$ with the following two properties: 
\beginitems
\itemrm{(i)}
there is a constant $m$ such that
for every $x \in \vertex(G)$, 
$$
\big|\big\{i \st x \in \vertex(H_i)\big\}\big| = m
$$
\enditems
and
\beginitems
\itemrm{(ii)} 
for every $e \in \edges(G)$, 
$$
w(e)
\ge
{1 \over m} \sum_{i \st e \in \edges(H_i)} w_i(e)
\,.
$$
\enditems
Then for all $t > 0$, we have 
$$
{1 \over |G|} \sum_{x \in \vertex(G)} p_t(x; G)
\le
{1 \over \sum_{j=1}^k |H_j|}
\sum_{i=1}^k  \sum_{x \in \vertex(H_i)} p_t(x; H_i)
\,.
$$
\endprocl

We will use the notation $A \le B$ for self-adjoint operators $A$ and $B$
to mean that $B - A$ is positive semidefinite.
Sometimes we regard the edges of a graph as oriented, where we choose one
orientation (arbitrarily) for each edge. In particular, we do this whenever
we consider the $\ell^2$-space of the edge set of a graph.
In this case, we denote the tail
and the head of $e$ by $\etail e$ and $\ehead e$.
Define $d_G \colon \ell^2\big(\vertex(G)\big) \to \ell^2\big(\edges(G)\big)$ by
$$
d_G(a)(e)
:=
\sqrt{w(e)}\, \big[a(\etail e) - a(\ehead e)\big]
\,.
$$
Then $\Delta_G = d_G^* d_G$.
Let $\Tr$ denote normalized trace of a square matrix, i.e., the average of
the diagonal entries.
We use $\tr$ for the usual trace.

\proof
Let $n := |G|$ and $N := \sum_{j=1}^k |H_j| = nm$.
Write $V := \vertex(G)$.
Let 
$$
W := \bigcup_{i=1}^k \vertex(H_i) \times \{i\}
\,,
$$
so that $|W| = N$.
Suppose that $\Phi \colon \LL\big(\ell^2(W)\big) \to \LL\big(\ell^2(V)\big)$ 
is a positive unital linear map, i.e., a linear map
that takes positive operators to positive
operators and takes the identity map to the identity map.
(Here, a positive operator means positive semidefinite.)
Theorem 3.9 of \ref b.Ant:Jensen/
says that
$$
\Tr f\big(\Phi(A)\big)
\le
\Tr \Phi\big(f(A)\big)
\label e.jensen
$$
for self-adjoint operators $A \in \LL\big(\ell^2(W)\big)$ and functions $f
\colon \R \to \R$ that
are convex on the convex hull of the spectrum of $A$.
(In fact, those authors show the more general inequality 
$\Tr g\big(f\big(\Phi(A)\big)\big)
\le
\Tr g\big(\Phi\big(f(A)\big)\big)$ 
for every increasing convex $g$.)

We apply this as follows.
Write 
$$
\Gamma(x) := \{ i \st x \in H_i \}
$$
and define $\phi \colon \ell^2(V) \to \ell^2(W)$ by linearity and the
requirement that 
$$
\phi(\II x)
:= 
{1 \over \sqrt m} \sum_{i \in \Gamma(x)} \II{(x, i)}
\,.
$$
Then $\phi^* \phi$ is the identity map by hypothesis (i).
Define $\Phi \colon \LL\big(\ell^2(W)\big) \to \LL\big(\ell^2(V)\big)$ by 
$
\Phi T
:=
\phi^* T  \phi
$.
Then $\Phi$ is a positive unital map.
Regard $H_i \times \{i\}$ as a graph with weights $w_i$ and a corresponding
Laplacian matrix $\Delta_i$. 
Consider the following map $A \in \LL\big(\ell^2(W)\big)$:
$$
A := \bigoplus_{i = 1}^k \Delta_i
\,.
$$
Hypothesis (ii) guarantees that 
$$
\Delta_G 
\ge 
\Phi(A)
\,.
\label e.bigger
$$
To see this, let $b \in \ell^2(V)$. 
We have 
$$
\bigpip{\Delta_G(b), b}
=
\|d_G b\|^2
\label e.1
$$
and 
$$
\bigpip{\Phi A(b), b}
=
\bigpip{\phi^* A \phi b, b}
=
\bigpip{A \phi b, \phi b}
\,.
\label e.2
$$
Write $b_i$ for the orthogonal projection of $\phi b$ to
$\ell^2\big(\vertex(H_i) \times \{i\}\big)$, so that
$\phi b = \sum_{i=1}^k b_i$.
Note that $b_i(x, i) = b(x)/\!\sqrt m$ for $x \in \vertex(H_i)$.
Thus, we have
$$\eqaln{
\bigpip{A \phi b, \phi b}
&=
\sum_{i=1}^k \bigpip{\Delta_i b_i, b_i}
=
\sum_{i=1}^k \|d_{H_i} b_i\|^2
=
{1 \over m}
\sum_{i=1}^k \sum_{e \in \edges(H_i)} {w_i(e)\over w(e)} \big|d_{G} b(e)\big|^2
\cr&\le
\sum_{e \in \edges(G)} \big|d_{G} b(e)\big|^2
=
\|d_G b\|^2
}$$
by hypothesis.
Combining this with \ref e.1/ and \ref e.2/, we get our claimed inequality
\ref e.bigger/.

Since \ref e.bigger/ implies, by the minimax principle,
that the eigenvalues of $\Delta_G$ are at
least the corresponding eigenvalues of $\Phi(A)$,
we have 
$$
\Tr f(\Delta_G)
\le
\Tr f\big(\Phi(A)\big)
$$
for every decreasing function $f$.
(We have strict inequality if $f$ is strictly decreasing and we have strict
inequality in \ref e.bigger/.)
Use $f(s) := e^{-t s}$ in this and in \ref e.jensen/ to obtain 
$$
\Tr f(\Delta_G)
\le
\Tr \Phi\big(f(A)\big)
\,.
\label e.done
$$
The left-hand side equals 
$$
{1 \over n} \sum_{x \in \vertex(G)} p_t(x; G)
\,.
$$
We claim that the right-hand side equals
$$
{1 \over N}
\sum_{i=1}^k  \sum_{x \in \vertex(H_i)} p_t(x; H_i)
\,,
$$
which will complete the proof of the theorem.

Another way to state our claim is that
$$
\Tr \Phi\big(f(A)\big)
=
{1 \over N} \sum_{i=1}^k \tr f(\Delta_i)
\,.
$$
Now, since $f(A) = \bigoplus_{i=1}^k f(\Delta_i)$, we have 
$$\eqaln{
\Tr \Phi\big(f(A)\big)
&=
{1 \over n} \sum_{x \in V}  \bigpip{f(A) \phi \II{x}, \phi \II{x}}
\cr&=
{1 \over N} \sum_{x \in V}  \sum_{i \in \Gamma(x)} \sum_{j \in \Gamma(x)}
\bigpip{f(A) \II{(x, i)}, \II{(x, j)}}
\cr&=
{1 \over N} \sum_{x \in V}  \sum_{i \in \Gamma(x)} \sum_{j \in \Gamma(x)}
\bigpip{f(\Delta_i) \II{(x, i)}, \II{(x, j)}}
\cr&=
{1 \over N} \sum_{x \in V}  \sum_{i \in \Gamma(x)} 
\bigpip{f(\Delta_i) \II{(x, i)}, \II{(x, i)}}
\cr&=
{1 \over N} \sum_{i=1}^k \tr f(\Delta_i)
\,.
\Qed
}$$

A similar proof shows that if $f$ is any decreasing convex function and $H$
fractionally tiles $G$, then 
$$
\Tr f(\Delta_G) \le \Tr f(\Delta_H)
\,.
\label e.eigs
$$
However, it is {\it not\/} true that this inequality holds whenever $G \dom
H$; a counter-example is provided by taking $f(t) := (4-t)^+$ and $G$, $H$
the graphs shown in \ref f.cex-dom/.
Possibly, however, it holds whenever $G \dom H$ and $H$ is transitive; this
is not hard to verify when $H$ is an edge.

\midinsert   
      \line{%
      \hbox to 3.25truein{\hfill\Size y1.8 \epsfbox{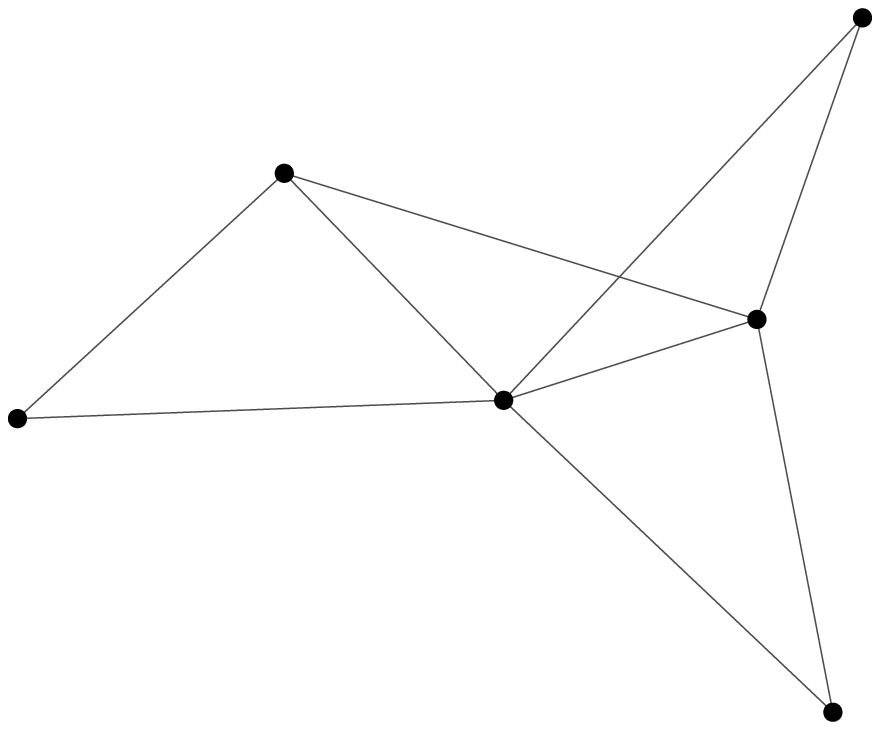}\hfill}%
      \hbox to 3.25truein{\hfill\Size y1.2 \raise
      .3truein\vbox{\epsfbox{sails2.eps}}\hfill}}%
\caption{\hsize=5truein%
\vtop{\noindent \figlabel{cex-dom}\enspace 
The graph $G$ on the left dominates the graph $H$ on the right.
}}
\endinsert

\procl r.sharp
\ref t.ftr/ is sharp in the following sense.
If the inequality in (ii) holds in the opposite direction for all edges
with strict inequality at least once, then the conclusion fails for all $t$
sufficiently close to 0. This is because both sides equal 1 for $t = 0$,
whereas the derivative of the left-hand side at $t = 0$ is 
$$
- \Tr \Delta_G
= -{2 \over |G|} \sum_{e \in G} w(e)
=
- {2 \sum_{e \in G} m w(e) \over \sum_i |H_i|}
$$
and the derivative of the right-hand side
at $t = 0$ is
$$
- {2 \sum_i \sum_{e \in H_i} w_i(e) \over \sum_i |H_i|}
\,.
$$
\endprocl

For special functions $f$, we can establish that domination is sufficient
for \ref e.eigs/.
A continuous function $f \colon (0, \infty) \to \R$ is called \dfn{operator
monotone on $(0, \infty)$} if for any bounded self-adjoint operators $A, B$
with spectrum in $(0, \infty)$ and $A \le B$, we have $f(A) \le f(B)$.
For example, \ref b.lowner/ proved that
the logarithm is an operator monotone function on $(0, \infty)$
\medbigp{see also Chapter V of \ref b.Bhatia/}.

\procl p.opmon
If $f$ is any operator monotone increasing function on $(0, \infty)$
and $G$ dominates $H$, then 
$$
\Tr f(\Delta_G + t) \ge \Tr f(\Delta_H + t)
\label e.opmon
$$
and
$$
\det (\Delta_G + t I)^{1/|G|}
\ge
\det (\Delta_H + t I)^{1/|H|}
\label e.charpoly
$$
for $t > 0$.
\endprocl

\proof
Fix $t > 0$ and define $g(s) := f(s+t)$.
Consider a copy $K$ of $H$ in $G$ and some vertex $x \in \verts(K)$. 
We have $\Delta_G \ge \Delta_K \oplus
\bfz$. Therefore $g(\Delta_G) \ge g(\Delta_K \oplus \bfz) = g(\Delta_K)
\oplus g(0) I$. Comparing the $(x, x)$-entries, we obtain $g(\Delta_G)(x,
x) \ge g(\Delta_K)(x, x)$.

The definition of $G \dom H$ is that there is
a certain coupling of copies of $(H, Y)$ and $(G, X)$ with $Y$
mapping to $X$; for each such copy $(K, X)$, apply the preceding
inequality and take expectation. This gives \ref e.opmon/.

Taking $f = \log$ yields \ref e.charpoly/. 
\Qed

\bsection{Fractional Tiling and Independent Sets}{s.fracind}

There are some easy results that follow from Shearer's inequality (\ref
B.CGFS/), which states the following:

\procl t.shearer 
Given discrete random variables $X_1$, \dots, $X_k$ and $S \subseteq [1, k]$,
write $X_S$ for the random variable $\Seq{X_i \st i \in S}$.
Let $\ev S$ be a collection of subsets of $[1, k]$ such that each integer
in $[1, k]$ appears in exactly $r$ of the sets in $\ev S$.
Then
$$
r \,\ent(X_1, \ldots, X_k)
\le
\sum_{S \in \ev S} \ent(X_S)
\,.
$$
\endprocl

Here, $\ent(X) := - \sum_x \P[X = x] \log \P[X = x]$
denotes the entropy of a discrete random variable $X$.

A set of vertices in a graph is \dfn{independent} if no pair in the set is
adjacent.
A \dfn{homomorphism} from $G$ to $H$ is a map from $\verts(G)$ to
$\verts(H)$ that sends adjacent vertices to adjacent vertices.
If $w \colon \verts(H) \to (0, \infty)$ is a weight function, then the
\dfn{weight} of a map $\phi \colon \verts(G) \to \verts(H)$ is $\prod_{x
\in \verts(G)} w\onebig(\phi(x)\onebig)$.
The total weight of a set of such maps is just the sum of the weights of
the individual maps.

\procl p.indepsets
Let $f(G)$ denote one of the following:
\beginbullets

the number of independent sets in $G$;

the number of proper colorings of $G$ with a fixed number of colors;

the total weight
of the homomorphisms of\/ $G$ to a fixed graph $F$ with arbitrary
positive weights on the vertices of $F$.

\endbullets
\noindent
If $H$ fractionally tiles $G$, then 
$$
f(G)^{1/|G|} 
\le
f(H)^{1/|H|} 
\,.
\label e.vxfrac
$$
\endprocl

\proof
For each of the things we count, the restriction of one of them in $G$ to a
copy of $H$ is also one of them for $H$.
Thus, \ref e.vxfrac/ is immediate from Shearer's inequality:
For example, if $S$ is a random uniform independent set in $G$, then its
entropy $h(S)$ equals $\log f(G)$.
Let $Z_x := \II{x \in S}$.
Let $H_j$ ($j \in J$) be the copies of $H$ that fractionally tile $G$.
Each vertex of $G$ belongs to exactly $|J| \cdot |H| /
|G|$ of these copies of $H$.
Since the restriction of $S$ to $\vertex(H)$ is an independent set in $H$,
we have 
$$\eqaln{
\log f(G)
&=
h\big(\Seqbig{Z_x \st x \in \vertex(G)}\big)
\cr&\le 
{|G| \over |J| \cdot |H|} 
\sum_j h\big(\Seqbig{Z_x \st x \in \vertex(H_j)}\big)
\cr&\le
{|G| \over |J| \cdot |H|} 
\sum_j \log f(H_j)
=
{|G| \over |H|} 
\log f(H)
}$$
by Shearer's inequality.

For weighted homomorphisms, standard techniques apply: it suffices to prove it
for rational weights or, by homogeneity, for integral weights.
Then we blow up each vertex of $F$ a
certain number of times to get an equivalent inequality for 
an unweighted graph, $F'$, which follows as above.
\medbigp{Here, $F'$ has vertex set $\big\{(x, i) \st x \in \verts(F),\, 
1 \le i \le w(x)\big\}$ and an
edge from $(x, i)$ to $(y, j)$ whenever $(x, y) \in \edges(F)$.}
\Qed

A similar inequality holds for fractional tilings by varied graphs, rather
than by a fixed graph.
That is, if $G$ is fractionally tiled by $H_1, \ldots, H_k$, meaning that
each $H_i$ is a subgraph of $G$ and each vertex of $G$ belongs to the same
number of $H_i$, then 
$$
f(G)^{1/|G|} 
\le
\Big(\prod_{j=1}^k f(H_j)\Big)^{1\big/\sum_{j=1}^k |H_j|}
\,.
$$
Of course, similar inequalities hold for hypergraphs.

We do not know when \ref e.vxfrac/ holds under the weaker assumption
that $G$ dominates $H$.
It does not hold when $f$ counts independent sets, as the
example of $G$ being a star and $H$ being an edge shows.
For proper colorings, however, it is easy to check that this inequality does
hold when $H$ is an edge.

The following is proved similarly to \ref p.indepsets/, but with basic
random variables representing edges rather than vertices.
In this proposition, we say that $H$ \dfn{fractionally edge-tiles} $G$ if
there is a set of copies of $H$ in $G$ such that each edge of $G$ belongs
to the same number of copies of $H$ in the set.

\procl p.edgesets
Let $f(G)$ denote one of the following:
\beginbullets

the number of acyclic orientations of $G$
(this is an evaluation of the Tutte polynomial, $T_G(2, 0)$);

the number of forests in $G$
(this is an evaluation of the Tutte polynomial, $T_G(2, 1)$);

the number of matchings in $G$.
\endbullets
\noindent
If $H$ fractionally edge-tiles $G$, then 
$$
f(G)^{1/|\edges(G)|} 
\le
f(H)^{1/|\edges(H)|} 
\,.
\label e.edgefrac
$$
\endprocl

We do not know when the inequality opposite to
\ref e.vxfrac/ holds under the weaker assumption
that $G$ dominates $H$ for the functions $f$ of \ref p.edgesets/.
It does not hold when $f$ counts matchings, as the
example of $G$ being a star and $H$ being an edge shows.
It might be the case that for any $x, y \ge 1$, we have $T_G(x,
y)^{1/|G|} \ge T_H(x, y)^{1/|H|}$ when $G$ dominates
$H$, or even that all the coefficients of $T_G(x+1, y+1)^{|H|} -
T_H(x+1, y+1)^{|G|}$ are non-negative.
Random testing of pairs $G \dom H$ supports the possibility that
all coefficients are non-negative in this difference.
Of course, such an inequality would imply \ref e.spnineq/.
When $H$ is a tree,
it is easy to prove this inequality.

We close with a few questions involving fractional tiling.

Let $f(G)$ be the number of matchings of $G$. Is
$$
f(G)^{1/|G|} 
\ge
f(H)^{1/|H|} 
\label e.match
$$
when $H$ fractionally tiles $G$?
After discussions with \'Ad\'am Tim\'ar, he proved this holds when $H$ is
an edge. To see this,
consider a maximal matching, $M$, of $G$. Let $W := \verts(G)
\setminus \verts(M)$. By definition, $W$ is an independent set.
In a fractional tiling of $G$ by $H$, let $K$ be the list of
all the copies of $H$ that use a vertex of $W$. This may include
repetitions. The fact that $H$ fractionally tiles $G$ combined with
Hall's marriage theorem allows us to conclude that there is a matching $M'$
that contains $W$ with $M'$ using only edges from $K$, and where each edge
in $M'$ intersects $W$. 
Now $M \cup M'$ is a subgraph of $G$ that satisfies the
inequality, as can be seen by considering the connected components of $M
\cup M'$.
That is, 
$$
f(G)^{1/|G|} 
\ge
f(M \cup M')^{1/|G|} 
\ge
2^{1/2}
=
f(H)^{1/|H|} 
\,.
$$

More generally, call a set of subgraphs of $G$ a \dfn{packing} if the
subgraphs are disjoint.
Let $f(G)$ be the number of packings of $G$ by copies of a
fixed graph $K$ (so when $K$ is an edge, this is the number of matchings).
Does \ref e.match/ hold when $H$ fractionally tiles $G$?
What about the particular case $K = H$?

\medbreak
\noindent {\bf Acknowledgements.}\enspace 
I am grateful to Oded Schramm for discussions at the start of this project
and to Prasad Tetali for conversations about \ref s.fracind/.
I thank Jeff Kahn and \'Ad\'am Tim\'ar for permission to include their
proofs.

\def\noop#1{\relax}
\input \jobname.bbl

\filbreak
\begingroup
\eightpoint\sc
\parindent=0pt\baselineskip=10pt

Department of Mathematics,
831 E. 3rd St.,
Indiana University,
Bloomington, IN 47405-7106
\href{mailto:rdlyons@indiana.edu}{{\tt rdlyons@indiana.edu}}
\par
\href
{http://pages.iu.edu/\string~rdlyons/}
{{\tt http://pages.iu.edu/\string~rdlyons/}}

\endgroup

\bye